\documentclass[twoside]{article}
\include{graphicx}
\title{The Stokes phenomenon associated with the periodic zeta function $F(a,s)$}
\author{R. B. Paris\\
\\
School of Engineering, Computing and Applied Mathematics,\\
 University of Abertay Dundee, Dundee DD1 1HG, UK\\
E-Mail: r.paris@abertay.ac.uk}
\begin{document}
\def\f#1#2{\mbox{${\textstyle \frac{#1}{#2}}$}}
\def\dfrac#1#2{\displaystyle{\frac{#1}{#2}}}
\def\boldal{\mbox{\boldmath $\alpha$}}
\newcommand{\bee}{\begin{equation}}
\newcommand{\ee}{\end{equation}}
\newcommand{\lam}{\lambda}
\newcommand{\ka}{\kappa}
\newcommand{\al}{\alpha}
\newcommand{\om}{\omega}
\newcommand{\Om}{\Omega}
\newcommand{\fr}{\frac{1}{2}}
\newcommand{\fs}{\f{1}{2}}
\newcommand{\g}{\Gamma}
\newcommand{\br}{\biggr}
\newcommand{\bl}{\biggl}
\newcommand{\ra}{\rightarrow}
\renewcommand{\topfraction}{0.9}
\renewcommand{\bottomfraction}{0.9}
\renewcommand{\textfraction}{0.05}
\newcommand{\mcol}{\multicolumn}
\date{}
\maketitle
\begin{abstract}
The exponentially improved large-$a$ expansion for the Hurwitz zeta function $\zeta(s,a)$ is exploited to
examine the expansion of the periodic zeta function $F(a,s)$ in the upper half-plane of the variable $a$. It
is shown that a double Stokes phenomenon takes place in the vicinity of the positive imaginary $a$-axis as $|a|\ra\infty$. This is a consequence of the fact that constituent parts of $F(a,s)$ involve two Hurwitz zeta functions resulting in two parallel Stokes lines at unit distance apart. Numerical calculations confirm the theoretical predictions.

\vspace{0.4cm}


\vspace{0.3cm}

\end{abstract}
\vspace{0.3cm}

\begin{center}
{\bf 1. \  Introduction}
\end{center}
\setcounter{section}{1}
\setcounter{equation}{0}
\renewcommand{\theequation}{\arabic{section}.\arabic{equation}}
It is a well-established result of asymptotic analysis that the Stokes phenomenon is associated with the smooth appearance of an exponentially small term in compound asymptotic expansions across certain rays (known as Stokes lines) in the complex plane. For a wide class of functions, particularly those satisfying second-order ordinary differential equations, the functional form of the coefficient multiplying such a subdominant exponential (a Stokes multiplier) is found to possess a universal structure represented to leading order by an error function,
whose argument is an appropriate variable describing the transition across the Stokes line \cite{B1}.

A function not satisfying a differential equation and which does not share this simple property is the logarithm of the gamma function. In \cite{PW}, Paris and Wood obtained the exponentially improved expansion of $\log\,\g(z)$ and showed that it involved not one but an {\it infinite number} of subdominant exponentials $e^{2\pi ikz}$ ($k=\pm1, \pm2, \ldots$). These exponentials are maximally subdominant on the Stokes lines $\arg\,z=\pm\fs\pi$, respectively, and steadily grow in magnitude in $|\arg\,z|>\fs\pi$ to eventually combine to generate the poles of $\g(z)$ on the negative $z$-axis. These authors demonstrated that the Stokes multipliers associated with the leading exponentials 
(corresponding to $k=\pm 1$) undergo a smooth transition approximately described by an error function in the neighbourhood of $\arg\,z=\pm\fs\pi$.
Subsequently, Berry \cite{B2} showed, by a sequence of increasingly delicate subtractions of optimally truncated asymptotic series, that all the subdominant exponentials switch on smoothly across the Stokes lines in a similar manner; see also \cite[\S 6.4]{PK}.

An analogous refinement in the large-$a$ asymptotics of the Hurwitz zeta function $\zeta(s,a)$ was considered by the author in \cite{P}. Across the Stokes lines $\arg\,a=\pm\fs\pi$, there is a similar appearance of an infinite number of subdominant exponentials $e^{2\pi ika}$ ($k=\pm 1, \pm 2, \ldots $), each exponential being associated with its own Stokes multiplier. For large $|a|$, the Stokes multipliers associated with these exponentials   
undergo a smooth, but rapid, transition given approximately by
\bee\label{e11}
\hspace{4cm}\fs\pm\fs\mbox{erf}\,[(\theta\mp\fs\pi)\sqrt{\pi k|a|}],\qquad \theta=\arg\,a \quad(k=1, 2, \ldots )
\ee
in the neighbourhood of $\arg\,a=\pm\fs\pi$, respectively.

In this paper we consider the periodic zeta function $F(a,s)$ defined by \cite[\S 25.13]{DLMF}
\bee\label{e12}
F(a,s)=\sum_{n=1}^\infty \frac{e^{2\pi ia}}{n^s}\qquad (\Re (s)>0,\ 0<a<1;\  \Re (s)>1,\ a\in {\bf N}),
\ee
which equals $\zeta(s)$ when $a$ is an integer. Here, we allow $a$ to take on complex values in the upper half-plane $0<\arg\,a<\pi$, when the sum in (\ref{e12}) converges for all finite values of the parameter $s$.
The function $F(a,s)$ can be expressed in terms of the Hurwitz zeta function and, accordingly, its exponentially improved large-$a$ expansion also consists of an infinite number of subdominant exponentials $e^{2\pi ika}$ ($k=\pm 1, \pm 2,\ldots$). In the neighbourhood of the positive imaginary axis, it is found that the exponentials with $k\geq 1$ undergo a double Stokes phenomenon, since constituent parts of $F(a,s)$ are associated with two parallel Stokes lines at unit distance apart.
In Section 2, we summarise the exponentially improved expansion of $\zeta(s,a)$ obtained in \cite{P} and apply this in Section 3 to $F(a,s)$.

\vspace{0.6cm}

\begin{center}
{\bf 2. \ The exponentially improved expansion for $\zeta(s,a)$}
\end{center}
\setcounter{section}{2}
\setcounter{equation}{0}
\renewcommand{\theequation}{\arabic{section}.\arabic{equation}}
We present here a summary of the details of the exponentially improved expansion of the Hurwitz zeta function $\zeta(s,a)$ derived in \cite{P}. This function is defined by
\[\zeta(s,a)=\sum_{k=0}^\infty \frac{1}{(k+a)^s} \qquad (a\neq 0, -1, -2, \ldots )\]
when $\Re (s)>1$ and elsewhere by analytic continuation. It is convenient to introduce the function $Z(s,a)$ by
\bee\label{e21}
\zeta(s,a)=\frac{1}{2}a^{-s}+\frac{a^{1-s}}{s-1}+\frac{Z(s,a)}{\g(s)}.
\ee
Then it is well known that $Z(s,a)$ has the large-$a$ (Poincar\'e) expansion given by \cite[p.~610]{DLMF}
\bee\label{e22}
Z(s,a)\sim \sum_{k=1}^\infty\frac{B_{2k}}{(2k)!}\,\frac{\g(2k+s-1)}{a^{2k+s-1}}
\ee
valid as $|a|\ra\infty$ in $|\arg\,a|<\pi$,  where $B_{2k}$ denote the even-order Bernoulli numbers.
Since successsive $B_{2k}$ have opposite signs, all terms in (\ref{e22}) have the same phase on $\arg\,a=\pm\fs\pi$, with the result that the positive and negative imaginary $a$-axes are Stokes lines for $\zeta(s,a)$.

Let $\{N_k\}$ ($k\geq 1$) denote a set of positive integers, which for the moment are arbitrary but which will subsequently be chosen to be optimal truncation indices, and define $\nu_k:=2N_k+s$. Define the so-called terminant function $T_\nu(z)$ as a multiple of the incomplete gamma function $\g(a,z)$ by
\[T_\nu(z)=e^{\pi i\nu}\,\frac{\g(\nu)}{2\pi i}\,\g(1-\nu,z).\]
Then the exponentially improved large-$a$ expansion of $Z(s,a)$ is \cite[Eq.~(2.5)]{P}
\bee\label{e23}
Z(s,a)=(2\pi)^{2s}\sum_{k=1}^\infty k^{s-1}\left\{\frac{1}{\pi}\sum_{r=0}^{N_k-1}(-)^r\frac{\g(2r+s+1)}{(2\pi ka)^{2r+s+1}}+R_k(a;N_k)\right\}
\ee
for $|\arg\,a|<\pi$ and\footnote{When $s=-n$ ($n=1, 2, \ldots$) we have $\zeta(-n,a)=-B_{n+1}(a)/(n+1)$, where $B_n(a)$ is the Bernoulli polynomial.} for all $s$ satisfying $s\neq -1, -2, \ldots \,$, where
\bee\label{e24}
R_k(a;N_k)=e^{-\pi is}\{e^{2\pi ika} e^{\fr\pi is}T_{\nu_k}(2\pi ika)-e^{-2\pi ika} e^{-\fr\pi is}T_{\nu_k}(-2\pi ika)\}.
\ee

An important feature of (\ref{e23}) is that the standard Poincar\'e expansion in (\ref{e22}) has been decomposed into a $k$-sequence of component asymptotic series with scale $2\pi ka$, each associated with its own {\it arbitrary} truncation index $N_k$ and remainder term $R_k(a;N_k)$. From the large argument asymptotics of the incomplete gamma function $\g(a,x)\sim x^{a-1}e^{-x}$ as $|x|\ra\infty$ in $|\arg\,x|<\f{3}{2}\pi$ \cite[p.~179]{DLMF}, the sum over $k$ involving $R_k(a;N_k)$ is seen to be absolutely convergent, since the behaviour of the late terms is controlled by $k^{-1-2N_k}$. It therefore follows that the result (\ref{e23}) is exact and no further expansion process is required. 

We remark in passing that in the particular case when all the indices $N_k$ are chosen equal to $N$ the sum over $k$ in (\ref{e23}) can be evaluated as $\zeta(2r+2)$. Use of the identity $B_{2r}=2(-)^{r-1}
\zeta(2r) (2r)!/(2\pi)^{2r}$ and some straightforward manipulation then shows that the algebraic series reduces to a truncated form of the Poincar\'e expansion in (\ref{e22}) to yield
\bee\label{e25}
Z(s,a)=\sum_{r=1}^N \frac{B_{2r}}{(2r)!}\,\frac{\g(2r+s-1)}{a^{2r+s-1}}+(2\pi)^s\sum_{k=1}^\infty k^{s-1} R_k(a;N)
\ee
with the index of the terminant functions contained in $R_k(a;N)$ given by $\nu_k=2N+s$. The result (\ref{e25}) has been given previously in \cite[Eq. (21)]{KT}.

If the truncation indices $N_k$ are chosen to correspond to the optimal truncation values (i.e., truncation at or near the least term in the inner series over $r$ in (\ref{e23})) then it is easily shown that
\bee\label{e26}
N_k\simeq \pi k|a|.
\ee
In this case, $\nu_k=2\pi k|a|+O(1)$ and we see that the order and the argument of each terminant function appearing in $R_k(a;N_k)$ are approximately equal in the limit $|a|\ra\infty$. When $|\nu|\sim |z|\gg 1$, with $\phi=\arg\,z$, the function $T_\nu(z)$ possesses the asymptotic behaviour \cite{O}, \cite[\S 6.2.6]{PK}
\bee\label{e27}
T_\nu(z)\sim\left\{\begin{array}{lc}\!\! \dfrac{-ie^{(\pi-\phi)i\nu}}{1+e^{-i\phi}}\,\frac{e^{-z-|z|}}{\sqrt{2\pi |z|}}\{1+O(z^{-1})\}& -\pi+\epsilon\leq\phi\leq\pi-\epsilon \\
\\
\!\!\fs+\fs \mbox{erf}\,[c(\phi)(\fs|z|)^\fr]+O(z^{-\fr}e^{-\fr|z|c^2(\phi)}),& \epsilon\leq\phi\leq2\pi-\epsilon\end{array}\right.
\ee
where $\epsilon$ denotes an arbitrarily small positive quantity and $c(\phi)$ is defined implicitly by
\[\fs c^2(\phi)=1+i(\phi-\pi)-e^{i(\phi-\pi)}\]
with the branch for $c(\phi)$ chosen so that $c(\phi)\simeq\phi-\pi$ near $\phi=\pi$.

Application of (\ref{e27}) to the terminant functions in (\ref{e24}), as carried out in \cite{P}, then shows that the approximate functional form of the Stokes multipliers for large $|a|$ in the vicinity of $\arg\,a=\pm\fs\pi$ is given by the error function in (\ref{e11}).

\vspace{0.6cm}

\begin{center}
{\bf 3. \ The Stokes phenomenon associated with $F(a,s)$}
\end{center}
\setcounter{section}{3}
\setcounter{equation}{0}
\renewcommand{\theequation}{\arabic{section}.\arabic{equation}}
The periodic zeta function $F(a,s)$ defined in (\ref{e12}) can be expressed in terms of the Hurwitz zeta function and so consequently we can exploit the exponentially improved expansion of $\zeta(s,a)$ summarised in Section 2.
If, for convenience, we replace $s$ by $1-s$ then \cite[Eq.~(25.13.2)]{DLMF}
\bee\label{e31}
F(a,1-s)=\sum_{k=1}^\infty k^{s-1}e^{2\pi ika}=\frac{\g(s)}{(2\pi)^s}\{e^{\fr\pi is}\zeta(s,a)+e^{-\fr\pi is}\zeta(s,1-a)\}.
\ee
In this last reference this formula is stated for $0<a<1$ and $\Re\,(s)<0$, whereas in \cite[Eq.~(4.1)]{P} it is established for $0<\arg\,a<\pi$ and arbitrary\footnote{When $s=0, -1, -2, \ldots$ the sum on the left-hand side of (\ref{e31}) can be expressed in terms of Clausen's function \cite[pp.~115--119]{SC}.}
 finite $s$ provided $s\neq 0, -1, -2, \ldots\,$. When $|a|>1$, the $-a$ in the second zeta function must be interpreted as $ae^{-\pi i}$ for $a$ in the upper half-plane.

To simplify the presentation we set $a'=1-a$ and define
\[{\tilde F}(a,s):=F(a,1-s)-\frac{\g(s)}{(2\pi)^s}\left\{e^{\fr\pi is}\left(\frac{1}{2}a^{-s}+\frac{a^{1-s}}{s-1}\right)+e^{-\fr\pi is}\left(\frac{1}{2}a'^{-s}+\frac{a'^{1-s}}{s-1}\right)\right\}.\]
Then, from (\ref{e21}) and (\ref{e31}) we obtain
\bee\label{e32}
{\tilde F}(a,s)=(2\pi)^{-s}\{e^{\fr\pi is}Z(s,a)+e^{-\fr\pi is}Z(s,a')\}.
\ee
If we now substitute the expansion (\ref{e23}) into (\ref{e32}) we find
\[{\tilde F}(a,s)=e^{\fr\pi is}\bl\{\frac{1}{\pi}\sum_{k=1}^\infty \frac{1}{k^{2r+2}}\sum_{r=0}^{N_k-1} A_r(a)+\sum_{k=1}^\infty k^{s-1} R_k(a;N_k)\br\}\]
\bee\label{e33}
\hspace{1.2cm}+e^{-\fr\pi is}\bl\{\frac{1}{\pi}\sum_{k=1}^\infty \frac{1}{k^{2r+2}}\sum_{r=0}^{N'_k-1} A_r(a')+\sum_{k=1}^\infty k^{s-1} R_k(a';N'_k)\br\}
\ee
where we have set
\bee\label{e33a}
A_r(a):=\frac{(-)^r\g(2r+s+1)}{(2\pi a)^{2r+s+1}}
\ee
and we have allowed for (possibly) different optimal truncation indices $N_k$ and $N'_k$ for the above algebraic asymptotic series involving the terms $A_r(a)$ and $A_r(a')$.

The expression in (\ref{e33}) can be written in the form
\bee\label{e34}
{\tilde F}(a,s)=\frac{e^{\fr\pi is}}{\pi}\sum_{k=1}^\infty\sum_{r=0}^{N_k-1}\frac{A_r(a)}{k^{2r+2}}
+\frac{e^{-\fr\pi is}}{\pi}\sum_{k=1}^\infty\sum_{r=0}^{N'_k-1}\frac{A_r(a')}{k^{2r+2}}+\sum_{k=1}^\infty k^{s-1}{\cal R}_k,
\ee
where, from (\ref{e24}),
\begin{eqnarray}
{\cal R}_k&=&e^{\fr\pi is}R_k(a;N_k)+e^{-\fr\pi is}R_k(a';N'_k)\nonumber\\
&=&e^{2\pi ika}T_{\nu_k}(2\pi ika)-e^{2\pi ika-\pi is}T_{\nu_k}(-2\pi ika)\nonumber\\
&& \hspace{2cm}+e^{-2\pi ika-\pi is}T_{\nu'_k}(2\pi ika')-e^{2\pi ika-2\pi is}T_{\nu'_k}(-2\pi ika')\label{e35}
\end{eqnarray}
with $\nu_k=2N_k+s$, $\nu'_k=2N'_k+s$ and we have used the fact that $e^{\pm 2\pi ika'}=e^{\mp2\pi ika}$.
In the neighbourhood of $\arg\,a=\fs\pi$, we note that the phase of the argument $-2\pi ika'$ of the final terminant function in (\ref{e35}) is approximately $-\pi$. 

Use of the connection formula \cite[p.~260]{PK}
\bee\label{e36b}
T_\nu(ze^{-\pi i})=e^{2\pi i\nu} (T_\nu(ze^{\pi i})-1)
\ee
then enables us to express the remainders ${\cal R}_k$ as
\[{\cal R}_k=e^{2\pi ika}\{T_{\nu_k}(2\pi ika)-T_{\nu'_k}(2\pi ika\xi)+1\}\hspace{3cm}\]
\bee\label{e36}
\hspace{4cm}+e^{-2\pi ika-\pi is}\{T_{\nu_k}(-2\pi ika)-T_{\nu'_k}(2\pi ika')\},
\ee
where
\bee\label{e37}
\xi:=1-\frac{1}{a},\qquad \delta:=\arg\,\xi.
\ee
Near $\arg\,a=\fs\pi$, the arguments of the first two terminant functions in (\ref{e36}) are approximately $+\pi$ (when $|a|\gg 1$), whereas those of the third and fourth terminant functions are approximately zero.

\vspace{0.6cm}

\begin{center}
{\bf 4. \ Numerical calculations }
\end{center}
\setcounter{section}{4}
\setcounter{equation}{0}
\renewcommand{\theequation}{\arabic{section}.\arabic{equation}}
In order to display numerically the smooth appearance of the $n$th subdominant exponential $e^{2\pi ina}$
in the vicinity of $\arg\,a=\fs\pi$ (at fixed $|a|$), it is necessary to `peel off' from ${\tilde F}(a,s)$ the larger subdominant exponentials corresponding to $1\leq k\leq n-1$ and all larger terms of the asymptotic series in (\ref{e34}). We follow the procedure described in \cite{B2}; see also \cite[\S 6.4.3]{PK} for an account of this process. 

By reversal of the order of summation it is seen that the first double sum in (\ref{e34}) can be rearranged\footnote{The modified double series involves the function $\zeta(2r+2,m)$ itself; however, its evaluation is straightforward for $m=1, 2, \ldots\,$.} as
\begin{eqnarray*}
\frac{1}{\pi}\sum_{k=1}^\infty \sum_{r=0}^{N_k-1} \frac{A_r(a)}{k^{2r+2}}&=& \frac{1}{\pi}\bl\{\sum_{r=0}^{N_1-1}\sum_{k=1}^\infty+\sum_{r=N_1}^{N_2-1}\sum_{k=2}^\infty+\sum_{r=N_2}^{N_3-1}\sum_{k=3}^\infty+\cdots\br\}\frac{A_r(a)}{k^{2r+2}}\\
&=&\frac{1}{\pi}\sum_{m=1}^\infty \sum_{r=N_{m-1}}^{N_m-1} A_r(a)\,\zeta(2r+2,m),
\end{eqnarray*}
where $N_0=0$, with a similar result for the other double sum involving $A_r(a')$. Then
\[{\tilde F}(a,s)=\frac{e^{\fr\pi is}}{\pi}\sum_{k=1}^\infty \sum_{r=N_{k-1}}^{N_k-1}A_r(a)\,\zeta(2r+2,k)\hspace{5cm}\]
\bee\label{e41}
\hspace{2cm}
+\frac{e^{-\fr\pi is}}{\pi}\sum_{k=1}^\infty\sum_{r=N'_{k-1}}^{N'_k-1}A_r(a')\,\zeta(2r+2,k)
+\sum_{k=1}^\infty k^{s-1}{\cal R}_k,
\ee
where ${\cal R}_k$ is defined in (\ref{e35}) and $A_r(a)$ in (\ref{e33a}).
This result is absolutely convergent and holds for arbitrary positive truncations $N_k$ and $N'_k$.

We now choose the integers $N_k$ and $N'_k$ to be optimal truncation indices according to (\ref{e26}). 
We define the $n$th Stokes multiplier $S_n(\theta)$ (with $\theta=\arg\,a$) associated with the exponential $e^{2\pi ina}$ by subtracting from ${\tilde F}(a,s)$ the larger subdominant exponentials $1\leq k\leq n-1$ and the asymptotic series 
corresponding to $1\leq k\leq n$, viz.
\[{\tilde F}(a,s)=\frac{e^{\fr\pi is}}{\pi}\sum_{k=1}^n \sum_{r=N_{k-1}}^{N_k-1}A_r(a)\,\zeta(2r+2,k)+\frac{e^{-\fr\pi is}}{\pi}\sum_{k=1}^n\sum_{r=N'_{k-1}}^{N'_k-1}A_r(a')\,\zeta(2r+2,k)\]
\[+\sum_{k=1}^{n-1} k^{s-1}{\cal R}_k+n^{s-1} e^{2\pi ina} S_n(\theta).\]
It then follows that
\[S_n(\theta)=\frac{e^{-2\pi ina}}{n^{s-1}}\bl\{{\tilde F}(a,s)-\frac{e^{\fr\pi is}}{\pi}\sum_{k=1}^n\sum_{N_{k-1}}^{N_k-1}A_r(a)\,\zeta(2r+2,k)\hspace{4cm}\]
\bee\label{e44}
\hspace{5cm}-\frac{e^{-\fr\pi is}}{\pi}\sum_{k=1}^n\sum_{N'_{m-1}}^{N'_m-1}A_r(a')\,\zeta(2r+2,k)-\sum_{k=1}^{n-1}k^{s-1}{\cal R}_k\br\}
\ee
for $n=1, 2, \ldots\ $. Note that the weighting factor $n^{s-1}$ has not been incorporated into the Stokes multiplier.

The Stokes multiplier for the first exponential $e^{2\pi ia}$ can be expressed alternatively following (\ref{e25}) in the form
\bee\label{e43}
S_1(\theta)=e^{-2\pi ia}\bl\{{\tilde F}(a,s)-\frac{e^{\fr\pi is}}{\pi}\sum_{r=1}^{N_1} \frac{B_{2r}}{(2r)!}\,\frac{\g(2r+s-1)}{a^{2r+s-1}}-
\frac{e^{-\fr\pi is}}{\pi}\sum_{r=1}^{N'_1}\frac{B_{2r}}{(2r)!}\,\frac{\g(2r+s-1)}{a'^{2r+s-1}} \br\}.
\ee

In Fig.~1 we present the real part\footnote{There is also a small imaginary part to $S_n(\theta)$ that is not shown.} of $S_n(\theta)$ for $n=1$ and $2$ when $a=|a|e^{i\theta}$ and for different $s$ as a function of $\theta$ in the vicinity of the positive imaginary $a$-axis. 
In the computation of $S_n(\theta)$ it is necessary to compute the terms ${\cal R}_k$ ($1\leq k\leq n-1$) by means of (\ref{e35}) and the double sums to the required exponential accuracy. The optimal truncation indices
$N_k$ and $N'_k$ were obtained by inspection of the terms in the algebraic expansions.
In addition, when computing the terminant functions appearing in ${\cal R}_k$ one must use the connection formula (\ref{e36b}) once the argument of $z$ in $T_\nu(z)$ has exceeded $\pi$, since {\it Mathematica} only computes the value of the incomplete gamma function in the principal sector $-\pi<\arg\,z\leq\pi$.

It is seen
from the asymptotics of $T_\nu(z)$ in (\ref{e27}) that in the neighbourhood of $\arg\,a=\fs\pi$ the dominant contribution to ${\cal R}_n$ arises from the terms involving the exponential $e^{2\pi ina}$,
since the terms involving $e^{-2\pi ina}$ are O($(n|a|)^{-1/2} e^{-2\pi n|a|})$.
The Stokes multiplier $S_n(\theta)$ of the factor $n^{s-1}e^{2\pi ina}$ contained in the remainder term in (\ref{e36}) then has the leading form given by
\begin{eqnarray}
S_n(\theta)&\simeq& T_{\nu_n}(2\pi ina)-T_{\nu'_n}(2\pi ina\xi)+1\nonumber\\
&\sim&1+\mbox{erf}\,[(\theta-\fs\pi)\sqrt{\pi n|a|}]-\mbox{erf}\,[(\theta+\delta-\fs\pi)\sqrt{\pi n|a'|}]\label{e42}
\end{eqnarray}
as $|a|\ra\infty$ for $n=1, 2, \ldots\,$, where $\xi$ and $\delta$ are defined in (\ref{e37}).
The computed values (shown by dots) are compared with the approximation (\ref{e42}) in Fig.~1. 
\begin{figure}[ht]
\begin{center}
{\small($a$)}\includegraphics[width=0.40\textwidth]{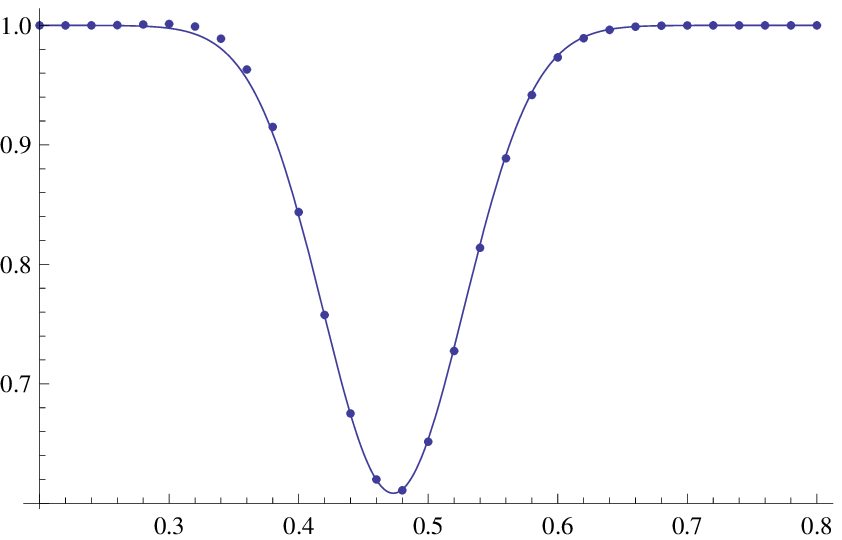}\hspace{1.2cm} {\small($b$)}\includegraphics[width=0.40\textwidth]{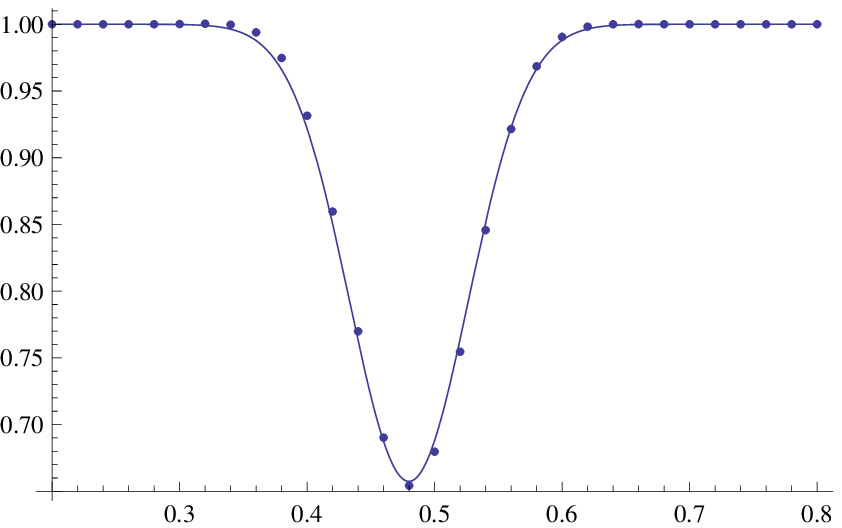}\vspace{0.8cm}
		
{\small($c$)}\includegraphics[width=0.40\textwidth]{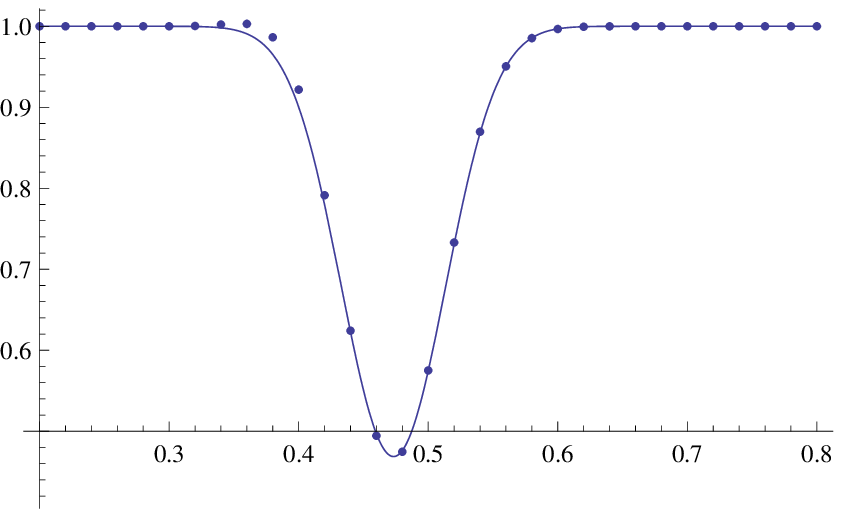}
\caption{\small{The variation of $\Re S_n(\theta)$ as a function of $\theta/\pi$ when 
	(a) $n=1$, $|a|=6$, $s=3$, with $N_1=N_1'=17$;  (b) $n=1$, $|a|=8$, $s=2+\fs i$, with $N_1=25$, $N_1'=24$;
	and (c) $n=2$, $|a|=6$, $s=2$, with $N_1=N_1'=18$, $N_2=36$, $N_2'=37$.
	The curves represent the approximate behaviour from (\ref{e42}) and the dots the exact values computed from (\ref{e44}) and (\ref{e43}).}}
\end{center}
\end{figure}
\vspace{0.6cm}

\begin{center}
{\bf 4. \ Concluding remarks }
\end{center}
\setcounter{section}{4}
\setcounter{equation}{0}
\renewcommand{\theequation}{\arabic{section}.\arabic{equation}}
The results in Fig.~1 show a good agreement between the computed values and the approximate functional form given by the right-hand side of (\ref{e42}). The transition of the exponentials in the vicinity of the positive imaginary axis occurs on an increasingly sharp scale (over a $\theta$-range proportional to $n^{-1/2}$). We remark that for the case $n=2$, the procedure has revealed the presence of an exponential of magnitude $e^{-24\pi}=O(10^{-33})$ hidden behind a larger exponential of magnitude $e^{-12\pi}$.

The most striking feature displayed in Fig.~1 is that the Stokes multipliers associated with ${\tilde F}(a,s)$ do not exhibit the simple switching on or off behaviour characteristic of the usual Stokes phenomenon, where the multiplier varies smoothly from approximately zero on one side of a Stokes line to approximately unity on the other side. This different behaviour associated with ${\tilde F}(a,s)$ is
produced by the fact that this function depends on two Hurwitz zeta functions, $\zeta(s,a)$ and $\zeta(s,1-a)$,
each with a different Stokes line in the upper half-plane: the first one on $\arg\,a=\fs\pi$ and the second on
$\Re\,(a)=1$, $\Im\,(a)>0$. As $\theta=\arg\,a$ is increased through positive values, the multiplier associated with $\zeta(s,1-a)$ switches off across $\Re\,(a)=1$ whereas the multiplier associated with $\zeta(s,a)$ switches on across $\Re\,(a)=0$. The result of these two conflicting transitions is that the multiplier for ${\tilde F}(a,s)$
undergoes a partial switch-off between these two Stokes lines, possessing the value approximately unity both before $\Re\,(a)=1$ and after $\Re\,(a)=0$, with a dip occurring in between these two lines. 

The minimum value of the approximate functional form of $S_n(\theta)$ in (\ref{e42}) depends only on $|a|$ and $n$. The location of the minimum at $\theta=\theta_0$ in the case $n=1$ is tabulated in Table 1 for different values of $|a|$ and also the corresponding minimum value $S_1(\theta_0)$.
It can be observed that the minimum $\theta_0$ is situated to the right of $\theta=\fs\pi$ and that $\theta_0\ra\fs\pi$ as $|a|\ra\infty$. The dip in $S_1(\theta)$ progressively decreases as $|a|\ra\infty$.
\begin{table}[th]
\begin{center}
\begin{tabular}{l|c|c||l|c|c}
\mcol{1}{c|}{$|a|$} & \mcol{1}{c|}{$\theta_0/\pi$} & \mcol{1}{c||}{$S_1(\theta_0)$} &
\mcol{1}{c|}{$|a|$} & \mcol{1}{c|}{$\theta_0/\pi$} & \mcol{1}{c}{$S_1(\theta_0)$} \\
[.1cm]\hline
&&&&\\[-0.25cm]
1 & $0.314363$ & $0.185967$ & 8 & $0.479894$ & $0.657472$\\
2 & $0.416139$ & $0.370072$ & 10& $0.483951$ & $0.691736$\\
4 & $0.459300$ & $0.529774$ & 15& $0.489331$ & $0.746192$\\
6 & $0.473089$ & $0.608463$ & 20& $0.492010$ & $0.779264$\\
[.1cm]\hline
\end{tabular}
\caption{\footnotesize{The location of the minimum $\theta_0$ of the leading approximation for $S_1(\theta)$ in (\ref{e42}) for different values of $|a|$.}}
\end{center}
\end{table}

\end{document}